\documentclass[10pt,leqno]{amsart}
\usepackage{eurosym}
\usepackage{indentfirst, amssymb,amsmath,amsthm,color}
\usepackage{csquotes}
\usepackage{hyperref}

\newtheorem{theo}{Theorem}[section]
\newtheorem{lem}{Lemma}[section]
\newtheorem{cor}{Corollary}[section]

\newcommand{\be}{\begin{equation}}
\newcommand{\ee}{\end{equation}}
\newcommand{\beas}{\begin{eqnarray*}}
\newcommand{\eeas}{\end{eqnarray*}}
\newcommand{\bea}{\begin{eqnarray}}
\newcommand{\eea}{\end{eqnarray}}

\numberwithin{equation}{section}
\input{tcilatex}

\begin{document}
\title[Uniqueness of an entire function with its shift difference ]{Some
results on uniqueness of an entire function with its shift difference}
\date{}
\author[Md. Adud, P. Saha and S. Tamang]{Md. Adud, Pratap Saha and Samten
Tamang }
\date{}
\address{Department of Mathematics, The University of Burdwan, Golapbag,
Burdwan, West Bengal, 713104, India}
\email{adud1993@gmail.com}
\address{Department of Mathematics, The University of Burdwan, Golapbag,
Burdwan, West Bengal, 713104, India}
\email{pratapsaha33@gmail.com}
\address{ Department of Mathematics, University of North Bengal, Raja
Rammohunpur, P.O.-N.B.U., Darjeeling, West Bengal, 734013, India}
\email{stamang@nbu.ac.in}
\maketitle

\begin{abstract}
In the paper, we have shown that every entire solution of the differential
difference equation $\Delta _{\eta }^{m}f-Q_{1}=(\Delta _{\eta
}^{k}f-Q_{2})e^{P}$ satisfy hyper order of $f=$degree of $P$ and using this
result we prove differential difference counter part Br\"{u}ck conjecture.
Also we have proved differential difference counter part Br\"{u}ck
conjecture for the special case of hyper order of $f<\frac{1}{2}.$
\end{abstract}

\let\thefootnote\relax\footnotetext{$^{1}$Corresponding Author} 
\footnotetext{%
} \footnotetext{%
} 
%%%%%%%%%%%%%%%%%%%%%%%%%%%%%%%%%%%%%%%%%%%%%%%%%%%%%%%%%%%%%%%%%%%%%%%%%%%%%%%%%%%%%%%%%%%%%%%%%%%%%%%%%%%%%%%

%%%%%%%%%%%%%%%%%%%%%%%%%%%%%%%%%%%%%%%%%%%%%%%%%%%%%%%%%%%%%%%%%%%%%%%%%%%%%%%%%%%%%%%%%%%%%%%%%%%%%%%%%%%%%%%%%%%%%%%%

\section{Introduction}

In this paper, we assume that the reader is familiar with the standard
results and theorems of Nevanlinna Theory, see, e.g., (\cite{HAYMAN}, \cite%
{HOLAND}, \cite{YANG AND YI}) . We also use some basis definitions,
notations and some results of the Wiman-Valiron theory (\cite{LAINE}). We
use $\sigma \left( f\right) $ and $\sigma _{2}\left( f\right) $ to denote
the order and the hyper order of an entire function $f,$ respectively where,%
\begin{equation*}
\sigma \left( f\right) =\underset{r\rightarrow \infty }{\lim \sup }\frac{%
\log N\left( r,f\right) }{\log r}\text{ and }\sigma _{2}\left( f\right) =%
\underset{r\rightarrow \infty }{\lim \sup }\frac{\log \log N\left(
r,f\right) }{\log r}.
\end{equation*}%
where $N\left( r,f\right) $ is the central index of $f.$

If $f\left( z\right) $ and $g\left( z\right) $ are two entire function and $%
a\in 
%TCIMACRO{\U{2102} }%
%BeginExpansion
\mathbb{C}
%EndExpansion
.$ We say $f$ and $g$ share the value $a$-CM(IM) if $f-a$ and $g-a$ have the
same zeros with counting(ignoring) multiplicities. In the uniqueness theory
Rubel and C.C. Yang (\cite{RUBEL}) proved that two entire function $f$ and $%
f^{^{\prime }}$ share two distinct value CM then $f$ and $f^{^{\prime }}$
are unique. Later Mues and Steinmetz(\cite{MUES}) improved the result and
prove that $f$ and $f^{^{\prime }}$ needs sharing two distinct value IM\ for
their uniqueness. In this case Br\"{u}ck made the following conjecture in
the year $1996.$

\textbf{Conjecture (\cite{BRUCK})}Let $f$ \ be an entire function such that
hyper order of $f$ is not a positive integer or $\infty $. If $f$ and $%
f^{^{\prime }}$ share one finite value $a$-CM, then%
\begin{equation*}
\frac{f^{^{\prime }}-a}{f-a}=c
\end{equation*}%
for some nonzero constant $c.$This conjecture was affirmed by Br\"{u}ck for $%
a=0$(\cite{BRUCK}) and by Gundersen and yang(\cite{Gundersen3}) for finite
order of $f$. This conjecture is still a open problem but many authors have
done a lot of work generalizing this conjecture. Recently many authors
started to consider the complex difference equation and shift difference for
the uniqueness of meromorphic functions. Let $f$ be a meromorphic function
and $\eta \in 
%TCIMACRO{\U{2102} }%
%BeginExpansion
\mathbb{C}
%EndExpansion
$, $m\in 
%TCIMACRO{\U{2115} }%
%BeginExpansion
\mathbb{N}
%EndExpansion
$. Denote $\Delta _{\eta }^{0}f=f(z),~\Delta _{\eta }f=\Delta _{\eta
}^{1}f=f(z+\eta )-f(z)$ and $\Delta _{\eta }^{m}f=\Delta _{\eta }\left(
\Delta _{\eta }^{m-1}f\right) =\overset{m}{\underset{j=0}{\sum }}\binom{m}{j}%
\left( -1\right) ^{m-j}f\left( m+j\eta \right) .$ Heittokanges et al proved
the following theorems.

\begin{theo}
\cite{HEITTOKANGAS}Let $f$ be a meromorphic function of order of growth%
\begin{equation*}
\sigma \left( f\right) =\underset{r\rightarrow \infty }{\lim \sup }\frac{%
\log T\left( r,f\right) }{\log r}<2,
\end{equation*}%
and let $\eta $ be a non-zero complex constant. If $f\left( z\right) $ and $%
f\left( z+\eta \right) $ share the values $a\in \mathbb{C}$ and $\infty $
CM, then%
\begin{equation*}
\frac{f\left( z+\eta \right) -a}{f\left( z\right) -a}=c,
\end{equation*}%
for some non-zero constant $c$.
\end{theo}

In (\cite{CHEN}), Chen proved the difference counterpart of the Br\"{u}ck
conjecture as follows:

\begin{theo}
\cite{CHEN}Let $f$ be a transcendental entire function of finite order that
is of a finite Borel exceptional value $\alpha \in \mathbb{C}$ and $\eta $
be a non-zero complex constant such that $f\left( z+\eta \right) \neq
\allowbreak f\left( z\right) .$ If $\Delta f\left( z\right) $ and $f\left(
z\right) $ share the value $a$ $\left( \neq \alpha \right) $CM, then%
\begin{equation*}
\frac{\Delta f\left( z\right) -a}{f\left( z\right) -a}=\frac{a}{a-\alpha }.
\end{equation*}
\end{theo}

Dong and Liu improve the above result and proved the following theorem

\begin{theo}
\cite{Dong and Liu}Let $f\left( z\right) $ be a transcendental entire
function of finite order and $\Psi \left( f\right) =\underset{j\in J}{\sum }%
a_{j}f^{\left( k_{j}\right) }\left( z+\eta _{j}\right) +\underset{i\in I}{%
\sum }b_{i}f\left( z+\xi _{i}\right) ,$ where $I,J$ are finite index sets, $%
a_{j},b_{i}$ are non zero constants, $k_{j}^{\prime }$s are positive
integers, $\eta _{j},\xi _{i}$ are complex constants such that $\underset{%
i\in I}{\sum }b_{i}=0.$ If $\Psi \left( f\right) $ and $f\left( z\right) $
share finite value $a-$CM and $f\left( z\right) $ has a finite exceptional
value $\alpha \left( \neq a\right) $ then%
\begin{equation*}
\left( i\right) ~\frac{\Psi \left( f\right) -a}{f\left( z\right) -a}=\text{%
constant,}
\end{equation*}%
if $a\neq 0$ and $\alpha $ is a Nevanlinna exceptional value%
\begin{equation*}
\left( ii\right) ~\frac{\Psi \left( f\right) -a}{f\left( z\right) -a}=\frac{a%
}{a-\alpha },
\end{equation*}%
if $\alpha $ is a borel exceptional value.
\end{theo}

In particular if we take $\Psi \left( f\right) =\left[ f\left( z+\eta
\right) -f\left( z\right) \right] ^{n}$ then we get the following corollary:

\begin{cor}
\cite{Dong and Liu}Let $f\left( z\right) $ be a transcendental entire
function of finite order that is of a finite Borel exceptional value and $%
f\left( z\right) \neq f\left( z+\eta \right) $ where $\eta $ is non-zero
complex constant. If $\Delta _{\eta }^{\left( n\right) }f=\left[ f\left(
z+\eta \right) -f\left( z\right) \right] ^{n}$ and $f\left( z\right) $ share
finite value $a\left( \neq \alpha \right) -$CM then%
\begin{equation*}
\frac{\Delta _{\eta }^{\left( n\right) }f-a}{f\left( z\right) -a}=\frac{a}{%
a-\alpha }.
\end{equation*}
\end{cor}

The aim of our paper is to confirm differential difference counter part of Br%
\"{u}ck conjecture for the case of hyper order of $f<\frac{1}{2}.$

\section{Main Results}

\begin{theo}
\label{theorem1}let $Q_{1}\left( z\right) $ and $\ Q_{2}\left( z\right) $
are two non-zero polynomials. If $f$ is a solution of the difference
differential equation 
\begin{equation}
\Delta _{\eta }^{m}f-Q_{1}=(\Delta _{\eta }^{k}f-Q_{2})e^{P}  \label{1}
\end{equation}%
where $P\left( z\right) $ is a polynomial of degree $n~$and $m,k$ are
non-negative integers such that $m-k>0$ then $\sigma _{2}\left( f\right) =$
degree of $P\left( z\right) .$
\end{theo}

\begin{theo}
\label{theorem2}let $Q_{1}\left( z\right) $ and $\ Q_{2}\left( z\right) $
are two non-zero polynomials and $P\left( z\right) $ be any polynomial. If $%
f $ is a non-constant solution of 
\begin{equation*}
\Delta _{\eta }^{m}f-Q_{1}=(\Delta _{\eta }^{k}f-Q_{2})e^{P}
\end{equation*}%
such that $\sigma _{2}\left( f\right) $ is not positive integer then, 
\begin{equation*}
\Delta _{\eta }^{m}f-Q_{1}=c(\Delta _{\eta }^{k}f-Q_{2})
\end{equation*}%
for some non zero constant $c.$
\end{theo}

\begin{cor}
$Q\left( z\right) $ and $\ P\left( z\right) $ are two non-zero polynomials
and if $f$ is a non constant solution of $\Delta _{\eta }f-Q=(f-Q)e^{P}$
such that $\sigma _{2}\left( f\right) $ is not positive integer then $\Delta
_{\eta }f-Q=c(f-Q),$ for some non zero constant $c.$
\end{cor}

%%%%%%%%%%%%%%%%%%%%%%%%%%%%%%%%%%%%%%%%%%%%%%%%%%%%%%%%%%%%%%%%%%%%%%%%%%%%%%%%%%%%%%%%%%%%%%%%%%%%%%%%%%%%%%%%%%%%%%

\begin{theo}
\label{theorem3}Let $f$ be a non constant entire function of hyper order $%
\sigma _{2}\left( f\right) <\frac{1}{2}$ and suppose $Q_{1}\left( z\right) $
and $\ Q_{2}\left( z\right) $ are two non-zero polynomials such that $\Delta
_{\eta }^{m}f-Q_{1}$ and $\Delta _{\eta }^{k}f-Q_{2}$ share $0-$CM, then 
\begin{equation*}
\Delta _{\eta }^{m}f-Q_{1}=c(\Delta _{\eta }^{k}f-Q_{2})
\end{equation*}%
for some non zero constant $c.$
\end{theo}

%%%%%%%%%%%%%%%%%%%%%%%%%%%%%%%%%%%%%%%%%%%%%%%%%%%%%%%%%%%%%%%%%%%%%%%%%%%%%%%%%%%%%%%%%%%%%%%%%%%%%%%%%%%%%%%%%%%%%%%%%%%%%%%%%%%%%%%%%%%%%
%%%%%%%%%%%%%%%%%%%%%%%%%%%%%%%%%%%%%%%%%%%%%%%%%%%%%%%%%%%%%%%%%%%%%%%%%%%%%%%%%%%%%%%%%%%%%%%%%%%%%%%%%%%%%%%%%%%%%%%%%%%%%%%%%%%%%%%%%%%%%

\section{Lemmas}

\begin{lem}
\label{lemma1}\cite{HOLAND}Let $P\left( z\right)
=p_{n}z^{n}+p_{n-1}z^{n-1}+...p_{1}z+p_{0,}~p_{n}\neq 0$ be a polynomial of
degree $n.$ Then there exist $R>0$ large enough such that for all $%
r=\left\vert z\right\vert >R$ the following inequalities hold%
\begin{equation*}
\left( 1-o\left( 1\right) \right) |p_{n}|r^{n}\leq \left\vert P\left(
z\right) \right\vert \leq \left( 1+o\left( 1\right) \right) |p_{n}|r^{n}.
\end{equation*}
\end{lem}

\begin{lem}
\label{lemma2}\cite{Ishizaki and Yanagihara}Let $f\left( z\right) $ be an
entire function then there exist a subset $E\subset \lbrack 1,\infty )$ with
finite logarithmic measure i.e., $\int_{E}\frac{dt}{t}<\infty $ and $%
r=\left\vert z\right\vert \notin \lbrack 0,1]\cup E$ with $\left\vert
f\left( z_{r}\right) \right\vert =M\left( r,f\right) $ such that%
\begin{equation*}
\frac{\Delta _{\eta }^{m}f}{f\left( z_{r}\right) }=\left( \frac{N\left(
r,f\right) }{z_{r}}\right) ^{m}\left( 1+o\left( 1\right) \right) ,~m\in 
%TCIMACRO{\U{2115} }%
%BeginExpansion
\mathbb{N}
%EndExpansion
,
\end{equation*}%
where $N\left( r,f\right) $ is the central index of $f.$
\end{lem}

\begin{lem}
\label{lemma3}\cite{barry}Let $f\left( z\right) $ be an entire function of
order $\sigma =\sigma \left( f\right) <\frac{1}{2},$ and let $A\left(
r\right) =\underset{\left\vert z\right\vert =r}{\inf }\log \left\vert
f\left( z\right) \right\vert $, $B\left( r\right) =\underset{\left\vert
z\right\vert =r}{\sup }\log \left\vert f\left( z\right) \right\vert .$ If $%
\sigma <\alpha <1,$ then%
\begin{equation*}
\underline{\log dens}\left\{ r:A\left( r\right) >\left( \cos \pi \alpha
\right) B\left( r\right) \right\} >1-\frac{\sigma }{\alpha },
\end{equation*}%
where%
\begin{equation*}
\underline{\log dens}\left( H\right) =\underset{r\rightarrow \infty }{\lim
\inf }\frac{\int_{1}^{r}\frac{\chi _{H}\left( t\right) }{t}dt}{\log r}
\end{equation*}%
and%
\begin{equation*}
\overline{\log dens}\left( H\right) =\underset{r\rightarrow \infty }{\lim
\sup }\frac{\int_{1}^{r}\frac{\chi _{H}\left( t\right) }{t}dt}{\log r},
\end{equation*}%
where $\chi _{H}\left( t\right) $ is the characteristic function of a set $%
H. $
\end{lem}

%%%%%%%%%%%%%%%%%%%%%%%%%%%%%%%%%%%%%%%%%%%%%%%%%%%%%%%%%%%%%%%%%%%%%%%%%%%%%%%%%%%%%%%%%%%%%%%%%%%%%%%%%%%%%%%%%%%%%%%%%%%%%%%%%%%%%
%%%%%%%%%%%%%%%%%%%%%%%%%%%%%%%%%%%%%%%%%%%%%%%%%%%%%%%%%%%%%%%%%%%%%%%%%%%%%%%%%%%%%%%%%%%%%%%%%%%%%%%%%%%%%%%%%%%%%%%%%%%%%

\section{Proof of the theorems}

\begin{proof}
(\textbf{Theorem \ref{theorem1}})Suppose $\ f$ is a polynomial. Then $\Delta
_{\eta }^{m}f-Q_{1}$ and $(\Delta _{\eta }^{k}f-Q_{2})$ are also polynomials
and so $P\left( z\right) $ must be a constants function. Therefore degree of 
$P\left( z\right) =0=\sigma _{2}\left( f\right) .$

Now let $f$ is a transcendental entire function.

\textbf{Case-I:} $\sigma \left( f\right) =\infty ,$ Then by Wiman-Valiron
theorem and Lemma(\ref{lemma2}) there exist a subset $E\subset \lbrack
1,\infty )$ with finite logarithmic measure such that $|z|=r\notin \lbrack
0,1]\cup E$ and $|f\left( z_{r}\right) |=m\left( r,f\right) ,$ we have%
\begin{equation}
\frac{\Delta _{\eta }^{m}f}{f\left( z_{r}\right) }=\left( \frac{N\left(
r,f\right) }{z_{r}}\right) ^{m}\left( 1+o\left( 1\right) \right) ,~m\in 
%TCIMACRO{\U{2115} }%
%BeginExpansion
\mathbb{N}
%EndExpansion
\label{2}
\end{equation}%
as \thinspace $|z_{r}|=r\rightarrow \infty ,$ and $N\left( r,f\right) $ is
the central index. Since $Q_{1}$ and $Q_{2}$ are non-zero polynomials,%
\begin{equation}
\underset{r\rightarrow \infty }{\lim }\frac{Q_{j}\left( z_{r}\right) }{%
f\left( z_{r}\right) }=\underset{r\rightarrow \infty }{\lim }\frac{%
Q_{j}\left( z_{r}\right) }{m\left( r,f\right) }=0  \label{3}
\end{equation}%
because $\sigma \left( f\right) =\infty ,$ $m\left( r,f\right) \rightarrow
\infty $ as $r\rightarrow \infty .~$Therefore from (\ref{1}), (\ref{2}) and (%
\ref{3}) we have 
\begin{eqnarray}
\frac{\frac{\Delta _{\eta }^{m}f}{f}-\frac{Q_{1}\left( f\right) }{f}}{\frac{%
\Delta _{\eta }^{k}f}{f}-\frac{Q_{2}\left( f\right) }{f}} &=&e^{P}  \notag \\
i.e.~\left( \frac{N\left( r,f\right) }{z_{r}}\right) ^{m-k}\left( 1+o\left(
1\right) \right) &=&e^{P}\text{ as }r\rightarrow \infty .  \label{4}
\end{eqnarray}%
Now since $P\left( z\right) $ is a polynomial of degree $n,$ let 
\begin{equation*}
P\left( z\right) =p_{n}z^{n}+p_{n-1}z^{n-1}+...p_{1}z+p_{0,}~p_{n}\neq 0
\end{equation*}%
Then by Lemma(\ref{lemma1}) for $|z|=r>R$ large enough we have 
\begin{equation}
\left( 1-o\left( 1\right) \right) |p_{n}|r^{n}\leq \left\vert P\left(
z\right) \right\vert \leq \left( 1+o\left( 1\right) \right) |p_{n}|r^{n}
\label{5}
\end{equation}%
Therefore by (\ref{4}) and (\ref{5}) we have%
\begin{eqnarray*}
\left( m-k\right) \log N\left( r,f\right) &\leq &\log \left\vert e^{P\left(
z_{r}\right) }\right\vert +\left( m-k\right) \log r+O\left( 1\right) \\
&\leq &\left\vert P\left( z_{r}\right) \right\vert +\left( m-k\right) \log
r+O\left( 1\right) \\
&\leq &\left( 1+o\left( 1\right) \right) \left\vert p_{n}\right\vert
r^{n}+\left( m-k\right) \log r+O\left( 1\right) .
\end{eqnarray*}%
Thus%
\begin{eqnarray}
\log \log N\left( r,f\right) &\leq &n\log r+\log \left\vert p_{n}\right\vert
+\log \log r+O\left( 1\right)  \notag \\
i.e.,~\sigma _{2}\left( f\right) &\leq &n=\text{ degree of }P\left( z\right)
\label{6}
\end{eqnarray}%
Now from (\ref{4}), taking principle log we have%
\begin{equation}
P\left( z_{r}\right) =\log \left\vert \left( \frac{N\left( r,f\right) }{z_{r}%
}\right) ^{k}\left( 1+o\left( 1\right) \right) \right\vert +i\arg \left\vert
\left( \frac{N\left( r,f\right) }{z_{r}}\right) ^{k}\left( 1+o\left(
1\right) \right) \right\vert  \label{7}
\end{equation}%
By (\ref{7}) and (\ref{5}) we have 
\begin{eqnarray*}
\left( 1-o\left( 1\right) \right) |p_{n}|r^{n} &\leq &P\left( z_{r}\right)
\leq \log \left\vert \left( \frac{N\left( r,f\right) }{z_{r}}\right)
^{m-k}\left( 1+o\left( 1\right) \right) \right\vert +O\left( 1\right) \\
&\leq &\left( m-k\right) \log N\left( r,f\right) +\left( m-k\right) \log
r+O\left( 1\right) \\
&\leq &2\left( m-k\right) \log N\left( r,f\right) .
\end{eqnarray*}%
Thus degree of $P\left( z\right) =n\leq \sigma _{2}\left( f\right) .$ Hence
Degree of $P\left( z\right) =\sigma _{2}\left( f\right) .$

\textbf{Case-II:} Suppose that $\sigma \left( f\right) <\infty .$ Then by
definition we have 
\begin{equation*}
\log N\left( r,f\right) =O\left( \log r\right)
\end{equation*}%
Now 
\begin{eqnarray*}
\left\vert P\left( z_{r}\right) \right\vert &=&\left\vert \log e^{P\left(
z_{r}\right) }\right\vert =\left\vert \left( m-k\right) \{\log N\left(
r,f\right) -\log r-i\theta \left( r\right) +O\left( 1\right) \}\right\vert \\
&\leq &O\left( \log r\right) ~\text{as }r\rightarrow \infty .
\end{eqnarray*}%
Since $P\left( z\right) $ is a polynomial and 
\begin{equation*}
\left\vert P\left( z_{r}\right) \right\vert \leq O\left( \log r\right)
\end{equation*}%
So $P\left( z\right) $ is constant. Hence degree of $P\left( z\right)
=0=\sigma _{2}\left( f\right) .$
\end{proof}

\begin{proof}
(\textbf{Theorem \ref{theorem2}})By theorem (\ref{theorem1}), $\sigma
_{2}\left( f\right) =$degree of $P\left( z\right) .$ Since According to our
hypothesis $\sigma _{2}\left( f\right) $ is not a positive integer.
Therefore our theorem follows.
\end{proof}

\begin{proof}
(\textbf{Theorem \ref{theorem3}})\bigskip Since $\Delta _{\eta }^{m}f-Q_{1}$
and $\Delta _{\eta }^{k}f-Q_{2}$ share $0-$CM therefore we can write 
\begin{equation}
\frac{\Delta _{\eta }^{m}f-Q_{1}}{\Delta _{\eta }^{k}f-Q_{2}}=e^{\beta },
\label{8}
\end{equation}%
Where $\beta $ is an entire function. If $f$ is a polynomial then $\Delta
_{\eta }^{m}f-Q_{1}$ and $\Delta _{\eta }^{k}f-Q_{2}$ are also polynomials.
Therefore $e^{\beta }=$constant, and our theorem follows. Now let $f$ is a
transcendantal entire function. We now consider the following two cases

\textbf{Case-I: }let $\sigma \left( f\right) <\infty .$ In this case the
proof can be carried out in the same way of case-II of Theorem(\ref{theorem1}%
).

\textbf{Case-II: }let $\sigma \left( f\right) =\infty .$ By the condition $%
\sigma _{2}\left( f\right) <\frac{1}{2}$ and from equation \ref{8} we can get%
\begin{equation}
\sigma \left( \beta \right) =\sigma _{2}\left( e^{\beta }\right) =\sigma
_{2}\left( \frac{\Delta _{\eta }^{m}f-Q_{1}}{\Delta _{\eta }^{k}f-Q_{2}}%
\right) \leq \sigma _{2}\left( f\right) <\frac{1}{2}.  \label{9}
\end{equation}%
By (\ref{9}) and Lemma(\ref{lemma3}) we can get a positive number $\alpha $
satisfying $\sigma <\alpha <\frac{1}{2}$ such that 
\begin{equation}
\underline{\log dens}\left\{ r:A\left( r\right) >\left( \cos \pi \alpha
\right) B\left( r\right) \right\} >1-\frac{\sigma }{\alpha },  \label{10}
\end{equation}%
where 
\begin{equation}
A\left( r\right) =\underset{\left\vert z\right\vert =r}{\inf }\log
\left\vert f\left( z\right) \right\vert  \label{11}
\end{equation}%
and%
\begin{equation}
B\left( r\right) =\underset{\left\vert z\right\vert =r}{\sup }\log
\left\vert f\left( z\right) \right\vert .  \label{12}
\end{equation}%
Now from (\ref{10}) and Wiman-Valiron theory we can get a subset $%
E_{1}\subset \left\{ r:A\left( r\right) >\left( \cos \pi \alpha \right)
B\left( r\right) \right\} \subset \lbrack 1,\infty )$ with finite
logarithmic measure such that equation (\ref{2}) holds for some $%
z_{r}=re^{i\theta \left( r\right) },~\left( \theta \left( r\right) \in
\lbrack 0,2\pi )\right) $ and $r\in S=:\left\{ r:A\left( r\right) >\left(
\cos \pi \alpha \right) B\left( r\right) \right\} \diagdown E_{1}$ and $%
M\left( r,f\right) =\left\vert f\left( z_{r}\right) \right\vert .$ Similarly
as in the proof of Theorem(\ref{theorem1}) equation (\ref{3}) and (\ref{4})
holds here. Now by equation (\ref{10}) we can write%
\begin{eqnarray}
A\left( z\right) &\leq &\log \left\vert \beta \left( z_{r}\right)
\right\vert \leq \left\vert \log \beta \left( z_{r}\right) \right\vert 
\notag \\
&=&\left\vert \log \log \frac{\Delta _{\eta }^{m}f-Q_{1}}{\Delta _{\eta
}^{k}f-Q_{2}}\right\vert  \notag \\
&=&\left\vert \log \log \frac{\frac{\Delta _{\eta }^{m}f}{f\left(
z_{r}\right) }-\frac{Q_{1}}{f\left( z_{r}\right) }}{\frac{\Delta _{\eta
}^{k}f}{f\left( z_{r}\right) }-\frac{Q_{2}}{f\left( z_{r}\right) }}%
\right\vert  \notag \\
&=&\left\vert \log \log \left( \left( \frac{N\left( r,f\right) }{z_{r}}%
\right) ^{m-k}\left( 1+o\left( 1\right) \right) \right) \right\vert  \notag
\\
&\leq &\log \log N\left( r,f\right) +O\left( 1\right) .  \label{13}
\end{eqnarray}%
As $r\in S$ and $r\rightarrow \infty .$\ Now from (\ref{12}) we can get 
\begin{equation}
\log M\left( r,\beta \right) =B\left( z\right) \leq \frac{A\left( z\right) }{%
\cos \pi \alpha }+O\left( 1\right) ~(r\in S).  \label{14}
\end{equation}%
Therefore by definition of hyper order and from equation (\ref{9}) we can
deduce $\log \log N\left( r,f\right) \leq O\left( \log r\right) ,$ as $r\in
S $ and $r\rightarrow \infty .$ Combining equation (\ref{13}) and (\ref{14})
we have 
\begin{equation*}
\log M\left( r,\beta \right) =O\left( \log r\right)
\end{equation*}%
as $r\in S$ and $r\rightarrow \infty .$ This imply that $\beta $ is a
polynomial. Thus from theorem (\ref{theorem1}) and equation \ref{8} and \ref%
{9} we have degree of $\beta =\sigma _{2}\left( f\right) <\frac{1}{2}.$
Which implies that $e^{\beta }$ is a constant so we get the conclusion of
our theorem.
\end{proof}

%%%%%%%%%%%%%%%%%%%%%%%%%%%%%%%%%%%%%%%%%%%%%%%%%%%%%%%%%%%%%%%%%%%%%%%%%%%%%%%%%%%%%%%

\section{Acknowledgements}

The authors are grateful to the anonymous referees for their valuable
suggestions which considerably improved the presentation of the paper.

For the first author this research is supported by University Grant
Commission(UGC), Grant No. 424608 (CSIR UGC NET 2017 JUNE), and the second
author is supported by the Council of Scientific and Industrial Research
(CSIR-HRDG), India, Grant No. 09/025(0255)/2018-EMR-I. 
%%%%%%%%%%%%%%%%%%%%%%%%%%%%%%%%%%%%%%%%%%%%%%%%%%%%%%%%%%%%%%%%%%%%%%%%%%%%%%%%%%%%%%%

\end{document}